\documentclass[11pt]{article}
  \RequirePackage{amsthm,amsmath}
     \usepackage{color,latexsym,amsfonts,amssymb}

\newcommand{\beqn}{\begin{eqnarray}}             
\newcommand{\eeqn}{\end{eqnarray}}               
\newcommand{\beq}{\begin{eqnarray*}}             
\newcommand{\eeq}{\end{eqnarray*}}

\newcommand{\ignore}[1]{}{}

\numberwithin{equation}{section}

\ignore{

\newbox\TempBox \newbox\TempBoxA

\newtheorem{lemma}{Lemma}[section]

\def\proclaim#1{\par \bigskip { #1}\bgroup\it\ }
\def\endproclaim{\egroup\par\bigskip}

\def\({\Big{(}}
\def\){\Big{)}}
\def\ep{\varepsilon}

\newtheorem{thm}{Theorem}[section]

\begin{document}

\bibliographystyle{agsm}

\usepackage[centertags]{amsmath}
\usepackage{amsfonts}
\usepackage{amssymb}
\usepackage{amsthm}
\usepackage{newlfont}

\titlepage
\textheight 235 truemm \textwidth 160 truemm \topmargin = -0.5cm
\oddsidemargin = 0.0cm \evensidemargin = -0.0cm

}

\newlength{\defbaselineskip}
\newcommand{\setlinespacing}[1]%
           {\setlength{\baselineskip}{#1 \defbaselineskip}}
\newcommand{\doublespacing}{\setlength{\baselineskip}%
                           {2.0 \defbaselineskip}}
\newcommand{\singlespacing}{\setlength{\baselineskip}{\defbaselineskip}}

\renewcommand{\v}{\vskip1cm}
\newcommand{\f}{\frac}
\newcommand{\wt}{\widetilde}
\renewcommand{\(}{\Big(}
\renewcommand{\)}{\Big)}
\renewcommand{\[}{\Big[}
\renewcommand{\]}{\Big]}
\renewcommand{\|}{\Big|}
\newcommand{\Xn}{X, X_1,\dots,X_n}
\newcommand{\Yn}{Y_1,\dots,Y_n}
\newcommand{\sn}{\sum_{i=1}^n}
\newcommand{\no}{\nonumber}
\newcommand{\ov}{\overline}
\newcommand{\la}{\label}
\newcommand{\be}{\begin{eqnarray}}
\newcommand{\ee}{\end{eqnarray}}
\newcommand{\bestar}{\begin{eqnarray*}}
\newcommand{\eestar}{\end{eqnarray*}}
\newcommand{\lam}{\lambda}
\newcommand{\Lam}{\Lambda}
\newcommand{\al}{\alpha}
\newcommand{\de}{\delta}
\newcommand{\De}{\Delta}
\newcommand{\ep}{\varepsilon}
\newcommand{\ga}{\gamma}
\newcommand{\Ga}{\Gamma}
\newcommand{\si}{\sigma}
\newcommand{\Si}{\Sigma}
\newcommand{\Th}{\Theta}
\newcommand{\pa}{\parallel}
\newcommand{\ka}{\kappa}
\newcommand{\nn}{\nonumber}

\allowdisplaybreaks

\def\Dto{\buildrel{D}\over \longrightarrow}

\newtheorem{thm}{Theorem}[section]
\newtheorem{lemma}{Lemma}[section]
\newtheorem{definition}{Definition}[section]
\newtheorem{cor}{Corollary}[section]
\newtheorem{prop}{Proposition}[section]
\newtheorem{rem}{Remark}[section]

\begin{document}

\title{Change in the mean in the domain of attraction of the normal law via Darling-Erd\H{o}s theorems}
\author{\bf Mikl\'{o}s Cs\"{o}rg\H{o}\thanks{Research supported by an NSERC Canada Discovery Grant at Carleton University.}\\
\footnotesize School of Mathematics and Statistics \\
\footnotesize Carleton University, 1125
Colonel By Drive, Ottawa, ON K1S 5B6, Canada\\
\bf mcsorgo@math.carleton.ca\\
\\
\bf Zhishui Hu\thanks{Partially supported by NSFC(No.10801122) and RFDP(No.200803581009), and by an NSERC Canada Discovery Grant of M.
Cs\"{o}rg\H{o} at Carleton University.}\\
 \footnotesize
Department of Statistics and Finance,  School of Management \\
\footnotesize  University of Science and Technology of China, Hefei,
Anhui 230026, China\\
\bf huzs@ustc.edu.cn}

\date{ }
\maketitle

\begin{center}
\begin{minipage}{130mm}
{{\bf Abstract.} This paper studies the problem of testing the null assumption of no-change in the mean of chronologically ordered independent
observations on a random variable $X$ {\it versus} the at most one change in the mean alternative hypothesis. The approach taken is via a
Darling-Erd\H{o}s type self-normalized maximal deviation between sample means before and sample means after possible times of a change in the
expected values of the observations of a random sample. Asymptotically, the thus formulated maximal deviations are shown to have a standard
Gumbel distribution under the null assumption of no change in the mean. A first such result is proved under the condition that $EX^2\log\log
(|X|+1)<\infty$, while in the case of a second one, $X$ is assumed to be in a specific class of the domain of attraction of the normal law,
possibly with infinite variance.}



\end{minipage}
\end{center}

\bigskip

\section{Introduction and main results}

Let $X, X_1,X_2,\cdots$ be non-degenerate independent identically distributed (i.i.d.) real-valued random variables (r.v.'s) with a finite mean
$EX=\mu$. We are interested in testing the null assumption
 \bestar H_0:~  X_1,X_2,\cdots X_n \mbox{~is a random sample on~} X \mbox{~with a
finite mean~} EX=\mu \eestar
 {\it versus} the ``at most one change
in the mean" (AMOC) alternative hypothesis \bestar && H_A:\mbox{~there is an integer~} k^*, 1\le k^*<n \mbox{~such that~}\\
&&~~~~~~EX_1=\cdots=EX_{k^*} \ne EX_{k^*+1}=\cdots =EX_n. \eestar

The hypothesized time $k^{*}$ of at most one change in the mean is usually unknown. Hence, given {\it chronologically ordered} independent
observables $X_1,X_2,\cdots, X_n, n\ge 1$, in order to test $H_0$ versus $H_A$, from a non-parametric point of view it appears to be reasonable
to compare the sample mean $(X_1+\cdots+X_k)/k =:S_k/k$ at any time $1\le k<n$ to the sample mean $(X_{k+1}+\cdots+X_n)/(n-k)
=:(S_n-S_k)/(n-k)$ after time $1\le k <n$ via functionals in $k$ of the family of the standardized statistics \be
 \Gamma_n(k)&:=&\Big(n\frac{k}{n}\Big(1-\frac{k}{n}\Big)\Big)^{1/2}\Big(\frac{S_k}{k}-\frac{S_n-S_k}{n-k}\Big)\nonumber\\
&=& \frac{1}{(\frac{k}{n}(1-\frac{k}{n}))^{1/2}}\Big(\frac{S_k}{n^{1/2}}-\frac{k}{n}\frac{S_n}{n^{1/2}}\Big),~~1\le k<n. \label{intr1}
 \ee
For instance, one would want to reject $H_0$ in favor of $H_A$ for large observed values of \be \Gamma_n:=\max_{1\le k< n} |\Gamma_n(k)|.
\label{intr2}\ee

 On the other hand, when assuming for example that the independent observables
 $X_1,\cdots, X_n, n\ge 1$, are $N(\mu, \sigma^2)$ random
variables, then we find ourselves modeling and testing for a parametric
 shift in the mean AMOC problem. It is, however, easy to check that,
when the variance $\sigma^2$ is known, then
 \be -2\log \Lambda_k=\frac{1}{\sigma^2}(\Gamma_n(k))^2,\label{intr3} \ee
where $\Lambda_k$ is the {\it likelihood ratio statistic} if the change in the mean occurs at $k^{*}=k$.
 Hence, the {\it maximally selected
likelihood ratio statistic} $\max_{1\le k<n}(-2\log \Lambda_k)$ will be large if and only if $\Gamma_n$ of (\ref{intr2}) is large. A similar
conclusion holds true if the variance $\sigma^2$ is an unknown but constant nuisance parameter (cf. Gombay and Horv\'{a}th (1994, 1996a,b), and
Cs\"{o}rg\H{o} and Horv\'{a}th (1997) [Section 1.4], and references therein). Namely in this case the maximally selected likelihood ratio
statistic $\max_{1\le k<n}(-2\log \Lambda_k)$ will be large if and only if \be \hat{\Gamma}_k:=\max_{1\le
k<n}\frac{1}{\hat{\sigma}_{k,n}}|\Gamma_n(k)|\label{intr4}\ee
 is large, where
\be \hat{\sigma}_{k,n}^2:=\frac{1}{n}\Big\{\sum_{1\le i\le k}\Big(X_i-\frac{S_k}{k}\Big)^2+\sum_{k< i\le
n}\Big(X_i-\frac{S_n-S_k}{n-k}\Big)^2\Big\}.\label{intr5}\ee

These conclusions, and further examples as well in Cs\"{o}rg\H{o} and Horv\'{a}th (1988) [Section 2], and in Cs\"{o}rg\H{o} and Horv\'{a}th
(1997) [Section 1.4] that are based on  Gombay and Horv\'{a}th (1994, 1996a,b), show that under the null hypothesis $H_0$ a large number of
parametric and nonparametric modeling of AMOC problems result in the same test statistic, namely that of (\ref{intr2}), or its variant in
(\ref{intr4}). Consequently, if the underlying distribution is not known, the just mentioned test statistics should continue to work just as
well when testing for $H_0$ versus $H_A$ as above. Furthermore, Brodsky and Darkhovsky (1993) argue quite convincingly in their Section 1.2
that detecting changes in the mean (mathematical expectation) of a random sequence constitutes one basic situation to which other changes in
distribution can be conveniently reduced. Thus $\Gamma_n$ and $\hat{\Gamma}_n$ gain a somewhat focal role in change-point analysis in general
as well. Studying the asymptotic behavior of these statistics is clearly of interest.

Let $S_0=0$, and for $n\ge 1$ define the sequence of tied-down partial sums processes
\be \label{intr6} Z_n(t):= \left\{\begin{array}{ll} (S_{[(n+1)t]}-[(n+1)t]S_n/n)/n^{1/2}, & 0\le t<1,\\
0, & t=1.
 \end{array}\right. \ee

In view of (\ref{intr1}), we are interested in exploring  the asymptotic behavior of the standardized sequence of stochastic processes \bestar
\Big\{\frac{1}{(t(1-t))^{1/2}}Z_n(t), 0\le t<1 \Big\}.\eestar We first note that \bestar
\sup_{0<t<1}\frac{1}{\sigma}|Z_n(t)|/(t(1-t))^{1/2}\eestar and, naturally, also the standardized statistics $\Gamma_n$ and $\hat{\Gamma}_n$
(cf. (\ref{intr2}) and (\ref{intr4})) converge in distribution to $\infty$ as $n\rightarrow \infty$  even if the null assumption of no change
in the mean is true. Hence, in order to secure nondegenerate limiting behavior under $H_0$, we seek appropriate renormalizations.

For example, it is proved in Cs\"{o}rg\H{o}, Szyszkowicz and Wang (2004) (cf. Corollary 5.2 in there) that, on {\it assuming $X$ to be in the
domain of attraction of the normal law (DAN), possibly with infinite variance}, then, as $n \rightarrow \infty$, \be \sup_{0<t<1}
\frac{1}{\hat{\sigma}_{[nt+1],n}} |Z_n(t)|/q(t) \stackrel{d}{\rightarrow} \sup_{0<t<1} |B(t)|/q(t),
 \label{intr7}
\ee where $\{B(t), 0\le t\le 1\}$ is a Brownian bridge, $\hat{\sigma}_{k,n},~1\le k\le n-1$ is as in (\ref{intr5}),
$\hat{\sigma}_{n,n}^2:=\frac{1}{n}\sum_{1\le i\le n}(X_i-\frac{S_n}{n})^2$,
\bestar  q(t):= \left\{\begin{array}{ll} (t\log\log (t^{-1}))^{1/2}, &  t\in (0,1/2],\\
((1-t)\log\log ((1-t)^{-1}))^{1/2}, & t\in [1/2, 1),
 \end{array}\right. \eestar
and $\log x:=\log (\max\{e,x\})$.

Large values of the statistics in (\ref{intr7}) indicate evidence against $H_0$. The weight function $q(\cdot)$ is to emphasize changes that
may have recurred near $0$ and $n$. We note in passing that the result in (\ref{intr7}) cannot be deduced via first proving a ``corresponding"
weak invariance principle on $D[0,1]$ (cf. Cs\"{o}rg\H{o} {\it et al.} (2004), Remark 5.2, as well as Corollaries 2 and 4 of  Cs\"{o}rg\H{o}
{\it et al.} (2008a) and their extension (46) in Theorem 4 of Cs\"{o}rg\H{o} {\it et al.} (2008b)). The applicability of (\ref{intr7}) is much
enhanced by Orasch and Pouliot (2004), tabulating functionals in weighted sup-norm.

An alternative way of studying change in the mean is via Darling-Erd\H{o}s type theorems. For example (cf. Theorems 2.1.2, A.4.2 and Corollary
2.1.2 in  Cs\"{o}rg\H{o} and Horv\'{a}th (1997)), {\it under $H_0$ with $EX^2\log\log (|X|+1)<\infty$, we have} \be \lim\limits_{n \rightarrow
\infty} P\Big(a(n)\max_{1\le k< n}\frac{1}{\hat{\sigma}_{k,n}}\Big(\frac{n^2}{k(n-k)}\Big)^{1/2}Z_n\Big(\frac{k}{n+1}\Big)\le
t+b(n)\Big)=\exp(-e^{-t}), ~~t \in \mathbb{R}, \label{intr9}\ee {\it where} \be
 a(n):=(2\log\log n)^{1/2}~\mbox{\it~and}~ ~b(n):=2\log\log n+\frac12
\log\log\log n-\frac12\log \pi.\label{intr10} \ee

In view of (\ref{intr7}), the aim of this paper is to explore the possibility of extending the result of (\ref{intr9}) to versions of
$Z_n(\frac{k}{n+1})$ under $H_0$ with $X\in $ DAN, for the sake of having an alternative approach to the sup-norm procedure of (\ref{intr7})
for studying the problem of a change in the mean in DAN, possibly with $EX^2=\infty$.

Define the family of statistics \be T_{k,n}=\frac{\frac{S_k}{k}-\frac{S_n-S_k}{n-k}} {\sqrt{\frac{\sum_{i=1}^k (X_i-S_k/k)^2}{k(k-1)}
+\frac{\sum_{i=k+1}^n (X_i-(S_n-S_k)/(n-k))^2}{(n-k)(n-k-1)}}},~~~2\le k\le n-2. \label{intr11} \ee

We note in passing that, on writing \be \tilde{\sigma}_{k,n}^2:=\frac{\sum_{1\le i\le k}\Big(X_i-\frac{S_k}{k}\Big)^2}{k(k-1)}+ \frac{\sum_{k<
i\le n}\Big(X_i-\frac{S_n-S_k}{n-k}\Big)^2}{(n-k)(n-k-1)},~~2\le k\le n-2, \label{addintr1}\ee we get \be T_{k,n}=
\frac{1}{\tilde{\sigma}_{k,n}}\Big(\frac{n}{k(n-k)}\Big)^{1/2}\Big(\frac{n^2}{k(n-k)}\Big)^{1/2}Z_n\Big(\frac{k}{n+1}\Big),~~2\le k\le n-2.
\label{addintr2}\ee

We note also that $(k(n-k)/n)\tilde{\sigma}_{k,n}^2$ is an unbiased estimator of $\sigma^2$ when $EX^2<\infty$.

Our first result is to say that, under the same moment condition for $X$, the self-normalized statistics $\max_{2\le k\le n-2} T_{k,n}$ behaves
 like $\max_{1\le k< n}\frac{1}{\hat{\sigma}_{k,n}}(\frac{n^2}{k(n-k)})^{1/2}Z_n(\frac{k}{n+1})$ does asymptotically
  (cf. our Theorem \ref{th2} and (\ref{intr9})). Our main result, Theorem \ref{th1}, however concludes the same asymptotic behavior for
  $\max_{1\le k<n} T_{k,n}$ for $X\in $
  DAN with possibly infinite variance.

\begin{thm} \label{th2}  Assume that $H_0$ holds and
\be EX^2\log\log (|X|+1)<\infty. \label{logcondition} \ee Then \bestar \lim\limits_{n \rightarrow \infty} P\Big(a(n)\max_{2\le k\le
n-2}T_{k,n}\le t+b(n)\Big)=\exp(-e^{-t}),~~t\in \mathbb{R}.\eestar
\end{thm}

Write $l(x):=E(X-\mu)^2I(|X-\mu|\le x)$. Assume that $X$ belongs to the domain of attraction of the normal law. Then $l(x)$ is a slowly varying
function as $x\rightarrow \infty$. Consequently, there exists some $a>1$ such that for any $x>a$ (see, for example, Galambos and Seneta
(1973)), \be \ell(x)=\exp\Big\{c(x)+\int_{a}^x \frac{\varepsilon(t)}{t}dt\Big\},\label{slowlyvarying} \ee where $c(x)\rightarrow c
(|c|<\infty)$ as $x\rightarrow \infty$ and $\varepsilon(t)\rightarrow 0$ as $t\rightarrow \infty$.

\begin{thm} \label{th1} Assume that $H_0$ holds and $l(x)$ is a slowly varying function at $\infty$ that, in terms of the representation
(\ref{slowlyvarying}), satisfies the additional conditions $c(x)\equiv c$ and $\varepsilon(t)\le C_0/\log t$ for some $C_0>0$, i.e., $X\in$
DAN, possibly with infinite variance, under the latter specific conditions on $l(x)$.
 Then,  for all $t \in \mathbb{R}$, \bestar \lim\limits_{n \rightarrow \infty} P\Big(a(n)\max_{2\le k\le n-2}T_{k,n}\le t+b(n)\Big)=\exp(-e^{-t}).\eestar
\end{thm}

{\bf Remark 1.} The additional conditions in Theorem \ref{th1} are satisfied by a large class of slowly varying functions, such as
$l(x)=(\log\log x)^{\alpha}$ and $l(x)=(\log x)^{\alpha}$, for example, for some $0<\alpha<\infty$.

{\bf Remark 2.}  Cs\"{o}rg\H{o},  Szyszkowicz and Wang (2003) obtained the follwoing Darling-Erd\H{o}s theorem for self-normalized sums: {\it
suppose that
 $H_0$ holds with $EX=0$ and $l(x)$ is a slowly varying function at $\infty$, satisfying
\be
 l(x^2)\le Cl(x) ~~\mbox{for some} ~~C>0. \label{lfun}
\ee Then, for every  $t \in \mathbb{R}$,}
 \bestar
\lim\limits_{n \rightarrow \infty} P\Big(a(n)\max_{1\le k\le n} S_k/V_k \le t+b(n)\Big)=\exp(-e^{-t}).
 \eestar
  If $l(x)$ has the representation
(\ref{slowlyvarying}) with $c(x)\equiv c$ and $\varepsilon(t)\le C_0/\log t$ for some $C_0>0$, then \bestar \frac{l(x^2)}{l(x)}
=\exp\Big\{\int_{x}^{x^2} \frac{\varepsilon(t)}{t}dt\Big \}\le \exp\Big\{C_0 \int_{x}^{x^2} \frac{1}{t\log t}dt\Big \}=2^{C_0}. \eestar So,
(\ref{lfun}) holds under the additional smoothness  conditions for $l(x)$ that are needed for  results like Lemma \ref{lemma2}, for example. On
the other hand, if $\varepsilon(x)=(\log x)^{-\alpha}$ for some $0<\alpha<1$, then $\lim\limits_{x\rightarrow \infty} l(x^2)/l(x)=\infty$,
i.e., (\ref{lfun}) fails. Thus, the additional conditions on $l(x)$ in Theorem \ref{th1} that are sufficient for having (\ref{lfun}), are seen
to be not far from being also necessary.

Before proving Theorems 1.1 and 1.2, we pose the following question.

\vskip 0.3cm
 {\bf Question 1.} In view of Theorems 1.1 and 1.2, one may like to know if the result of (\ref{intr9}) could also hold true when
replacing condition (\ref{logcondition}) by $X\in$ DAN, possibly with $EX^2=\infty$.

\vskip 0.3cm
 {\bf Question 2.} In view of having Theorems 1.1 and 1.2, one would hope to have (\ref{intr7}) in terms of $T_{k,n}$, i.e., when
replacing $\frac{1}{\hat{\sigma}_{[nt+1],n}}$ by $\frac{1}{\tilde{\sigma}_{[nt+1],n}}(\frac{n}{[nt+1](n-[nt])})^{1/2}$ on the left hand side of
(\ref{intr7}), with $\tilde{\sigma}_{k,n},~1\le k\le n-1$ defined as in (\ref{addintr1}) and $\tilde{\sigma}_{n,n}^2:=\frac{1}{n^2}\sum_{1\le
i\le n}(X_i-\frac{S_n}{n})^2$.

\vskip 0.3cm

As to these questions, it is clear from the respective proofs of (\ref{intr9}) (cf.
 Corollary 2.1.2 in  Cs\"{o}rg\H{o} and Horv\'{a}th (1997)) and Theorem \ref{th2} that,
 under the condition (\ref{logcondition}), the two estimators $\hat{\sigma}_{k,n}^2$
 and $({k(n-k)}/{n})\tilde{\sigma}_{k,n}^2$ of $\sigma^2$ are asymptotically equivalent.
 When $\mbox{Var} (X)=\infty$, this does not appear to be true any more, i.e., when these
 ``estimators" in hand are being used as self-normalizers. However, we could not resolve
 this problem as posed in the context of these two questions.

\section{Proofs of Theorems \ref{th2} and \ref{th1}}

Without loss of generality, in this section we assume that $\mu=0$.

\vskip 0.2cm
 \noindent {\bf Proof of Theorem \ref{th2}.} Write $K_n=\exp\{\log^{1/3} n\}$.  With $\tilde{\sigma}_{k,n}^2$ as in
(\ref{addintr1}), in view of (\ref{addintr2}),  at first, we prove that, as $n\rightarrow \infty$, \be \max_{K_n<
k<n-K_n}\Big|\frac{k(n-k)}{n}\tilde{\sigma}_{k,n}^2-\sigma^2\Big|=o_{P}((\log \log n)^{-1}).\label{hu1203}\ee Write $\tilde{b}_n=n/\log\log n$.
Then $\tilde{b}_n/n \downarrow$ and $\tilde{b}_n^2\sum_{i=n}^{\infty} \tilde{b}_i^{-2}=O(n)$. Noting that, for sufficiently large $n$, we have
\bestar
P\Big(|X^2-\sigma^2|>\tilde{b}_n\Big) &\le& P\Big(|X^2-\sigma^2|\log\log (|X^2-\sigma^2|+1)>\tilde{b}_n\log\log \tilde{b}_n\Big)\\
&\le& P\Big(|X^2-\sigma^2|\log\log (|X^2-\sigma^2|+1)>(1/2)n\Big),\eestar and $E|X^2-\sigma^2|\log\log (|X^2-\sigma^2|+1)<\infty$ (by the
assumption $EX^2\log\log (|X|+1)<\infty$), we conclude \bestar \sum_{n=1}^{\infty}P\Big(|X^2-\sigma^2|>\frac{n}{\log\log n}\Big)<\infty.
\eestar By Theorem 3 in Chow and Teicher (1978, page 126), we get
 \bestar \sum_{i=1}^k (X_i^2-\sigma^2)=o(k(\log\log k)^{-1})~~a.s.
~\mbox{~as~} k \rightarrow \infty. \eestar  Hence, by the classical Hartman-Wintner LIL,  as $k \rightarrow \infty$, we have
 \bestar \sum_{i=1}^k (X_i-S_k/k)^2-k\sigma^2=\sum_{i=1}^k (X_i^2-k\sigma^2)-S_k^2/k=o(k(\log\log
k)^{-1})~~a.s. \eestar Consequently,
 \bestar \max_{K_n< k \le n} \Big|\frac{1}{k}\sum_{i=1}^k
(X_i-S_k/k)^2-\sigma^2\Big|=o_P((\log\log n)^{-1}),\eestar and
 \bestar \max_{1\le k < n-K_n} \Big|\frac{1}{n-k}\sum_{i=k+1}^n
(X_i-(S_n-S_k)/(n-k))^2-\sigma^2\Big|=o_P((\log\log n)^{-1}).\eestar Hence (\ref{hu1203}) holds.

By Theorem 2.1.2 in Cs\"{o}rg\H{o} and Horv\'{a}th (1997), we have \bestar (2\log\log n)^{-1/2}\max_{1\le k<n} \Big(\frac{n}{k(n-k)}\Big)^{1/2}
\Big|S_k-\frac{k}{n}S_n\Big| \stackrel{P}{\rightarrow} \sigma. \eestar This, together with (\ref{hu1203}), implies \bestar
&&a(n)\Big|\max_{K_n<k<n-K_n} T_{k,n}-\frac{1}{\sigma}\max_{K_n<k<n-K_n} \Big(\frac{n}{k(n-k)}\Big)^{1/2}
\Big(S_k-\frac{k}{n}S_n\Big)\Big|\\
&&~~\le a(n)\max_{K_n<k<n-K_n} \Big(\frac{n}{k(n-k)}\Big)^{1/2} \Big|S_k-\frac{k}{n}S_n\Big|
\Big|\Big(\frac{k(n-k)}{n}\tilde{\sigma}_k^2\Big)^{-1/2}-\sigma^{-1}\Big|\\
&&~~ =o_P(1) (\log\log n)^{-1/2}\max_{1\le k<n} \Big(\frac{n}{k(n-k)}\Big)^{1/2} \Big|S_k-\frac{k}{n}S_n\Big| \stackrel{P}{\rightarrow}
0.\eestar
 Then  from the proof of Theorem A.4.2. in Cs\"{o}rg\H{o} and Horv\'{a}th (1997),  for all $t \in \mathbb{R}$,
it follows that
 \be \lim\limits_{n \rightarrow \infty} P\Big(a(n)\max_{K_n<k<n-K_n} {T}_{k,n}\le
t+b(n)\Big)=\exp(-e^{-t}). \label{addproof1}\ee

Similarly to the proof of (\ref{hu1201}) and (\ref{hu1202}) below, we get \be && a(n)\max_{2\le k\le K_n}
{T}_{k,n}-b(n)\stackrel{P}{\rightarrow} -\infty, \ee and \be && a(n)\max_{n- K_n\le k\le n-2} {T}_{k,n}-b(n)\stackrel{P}{\rightarrow}
-\infty.\label{addproof2}\ee

Now Theorem \ref{th2} follows from (\ref{addproof1})--(\ref{addproof2}). ~~$\Box$

\vskip 0.5cm

We continue with establishing three auxiliary lemmas for the proof of  Theorem \ref{th1}.

As in Cs\"{o}rg\H{o} {\it et al.} (2003), we start with putting $b=\inf\{x\ge 1;l(x)>0\}$ and \bestar \eta_n=\inf\Big\{s: s\ge
b+1,\frac{l(s)}{s^2}\le \frac{(\log\log n)^4}{n}\Big\}.\eestar Let \bestar
&&Z_j=X_jI(|X_j|>\eta_j),~~Y_j=X_jI(|X_j|\le \eta_j), ~~Y_j^{*}=Y_j-EY_j,\\
&&~~~~~S_n^{*}=\sum_{j=1}^n Y_j^{*}, ~~~~B_n^2=\sum_{j=1}^n EY_j^{*2},~~~~V_n^2=\sum_{j=1}^n X_j^2. \eestar Then, as $n\rightarrow \infty$,
$\eta_n \rightarrow \infty,~ nl(\eta_n)=\eta_n^2(\log\log n)^4$ for every large enough $n$ and $B_n^2\sim nl(\eta_n)$. As in Cs\"{o}rg\H{o}
{\it et al.} (2003), we may assume without loss of generality that \bestar B_n^2=nl(\eta_n)=\eta_n^2(\log\log n)^4 \mbox{~~for all~~} n\ge
1.\eestar

 Let $\{\tilde{X}, \tilde{X}_1,\tilde{X}_2,\cdots\}$  be a sequence of i.i.d. random variables  with $\tilde{X}\stackrel{d}{=} X$,
independently of $\{X, X_1, X_2, \cdots\}$. We  define $\tilde{S}_n,~\tilde{Z}_j, ~\tilde{Y}_j, ~\tilde{Y}_j^{*}, ~\tilde{S}_n^{*}$ and
$\tilde{V}_n$ similarly to ${S}_n,~{Z}_j, ~{Y}_j, ~{Y}_j^{*}, ~{S}_n^{*}$ and ${V}_n$.   Define \bestar
&& S_{k,n}=\left\{\begin{array}{ll}  {\frac{S_k}{k}-\frac{\tilde{S}_{n-[n/2]}+S_{[n/2]}-S_k}{n-k}}, & ~\,~\mbox{~~~~if~} 1\le  k\le n/2;\\
{\frac{S_{[n/2]}+\tilde{S}_{n-[n/2]}-\tilde{S}_{n-k}}{k}-\frac{\tilde{S}_{n-k}}{n-k}}, &  ~\,~\mbox{~~~~if~}  n/2<k<n, \end{array}\right. \\
&& S_{k,n}^*=\left\{\begin{array}{ll}  {\frac{S_k^*}{k}-\frac{\tilde{S}^*_{n-[n/2]}+S^*_{[n/2]}-S^*_k}{n-k}}, & ~\,~\mbox{~~~~if~}1\le k\le n/2;\\
{\frac{S^*_{[n/2]}+\tilde{S}^*_{n-[n/2]}-\tilde{S}^*_{n-k}}{k}-\frac{\tilde{S}^*_{n-k}}{n-k}}, &  ~\,~\mbox{~~~~if~}  n/2<k<n, \end{array}\right. \\
&& B_{k,n}^{2}=\left\{\begin{array}{ll}  \frac{B_k^2}{k^2}+\frac{B^2_{n-[n/2]}+B^2_{[n/2]}-B_k^2}{(n-k)^2}, & ~~\mbox{~if~} 1\le k\le n/2;\\
\frac{B^2_{[n/2]}+B^2_{n-[n/2]}-B^2_{n-k}}{k^2}+\frac{B^2_{n-k}}{(n-k)^2}, &  ~~\mbox{~if~} n/2<k<n,
\end{array}\right. \\
&& V_{k,n}^{2}=\left\{\begin{array}{l}  \frac{V_k^2}{k^2}- \frac{S_k^2}{k^3}+
\frac{\tilde{V}^2_{n-[n/2]}+V^2_{[n/2]}-V_k^2}{(n-k)^2}-\frac{(\tilde{S}_{n-[n/2]}+S_{[n/2]}-S_k)^2}{(n-k)^3},~~\mbox{~if~} 1\le k\le n/2;\\
\frac{V^2_{[n/2]}+\tilde{V}^2_{n-[n/2]}-\tilde{V}^2_{n-k}}{k^2}-
\frac{(S_{[n/2]}+\tilde{S}_{n-[n/2]}-\tilde{S}_{n-k})^2}{k^3}+\frac{\tilde{V}^2_{n-k}}{(n-k)^2}- \frac{\tilde{S}^2_{n-k}}{(n-k)^3},
\end{array}\right.  \\
&&~~~~~~~~~~~~~~~~~~~~~~~~~~~~~~~~~~~~~~~~~~~~~~~~~~~~~~~~~~~~~~~~~~~~~~~~~~~~~~~~\mbox{~~~if~} n/2<k<n.\\
  && \bar{V}_{k,n}^{2}=\left\{\begin{array}{l}  \frac{V_k^2}{k(k-1)}- \frac{S_k^2}{k^2(k-1)}+
\frac{\tilde{V}^2_{n-[n/2]}+V^2_{[n/2]}-V_k^2}{(n-k)(n-k-1)}-\frac{(\tilde{S}_{n-[n/2]}+S_{[n/2]}-S_k)^2}{(n-k)^2(n-k-1)},~~\mbox{~if~} 2\le k\le n/2;\\
\frac{V^2_{[n/2]}+\tilde{V}^2_{n-[n/2]}-\tilde{V}^2_{n-k}}{k(k-1)}-
\frac{(S_{[n/2]}+\tilde{S}_{n-[n/2]}-\tilde{S}_{n-k})^2}{k^2(k-1)}+\frac{\tilde{V}^2_{n-k}}{(n-k)(n-k-1)}-
\frac{\tilde{S}^2_{n-k}}{(n-k)^2(n-k-1)},
\end{array}\right.  \\
&&~~~~~~~~~~~~~~~~~~~~~~~~~~~~~~~~~~~~~~~~~~~~~~~~~~~~~~~~~~~~~~~~~~~~~~~~~~~~~~~~\mbox{~~~if~} n/2<k\le n-2.\eestar

  Clearly, with $\{T_{k,n}, 2\le
k\le n-2\}$ as in (\ref{intr11}), we have
 \bestar \{{T}_{k,n},  2\le
k\le n-2\}\stackrel{d}{=}\Big \{\frac{S_{k,n}}{\bar{V}_{k,n}},  2\le k\le n-2\Big\} ~~\mbox{for each}~~n\ge 4,\eestar where, and throughout,
$\stackrel{d}{=}$ stands for equality in distribution.

\begin{lemma} \label{lemma2} As $n \rightarrow \infty$, we have
\be &&\frac{l(\eta_n)-l(\eta_{n/(\log\log n)^5})}{l(\eta_n)}=o(1/\log\log n).\label{hu1204}\ee
\end{lemma}

\noindent {\bf Proof.}  Since \bestar && 1\ge \frac{l(\eta_{n/(\log\log n)^5})}{l(\eta_n)}\ge \exp\Big\{-C_0\int_{\eta_{n/(\log\log
n)^5}}^{\eta_n}\frac{1}{u\log
u}du\Big\}\\
&&~~\ge\exp\Big\{-C_0\frac{\eta_n}{\eta_{n/(\log\log n)^5}\log \eta_{n/(\log\log n)^5} }\Big\},\eestar and $\eta_n$ is a regularly varying
function with index $1/2$, for any $\varepsilon>0$, we have $\eta_n/ \eta_{n/(\log\log n)^5}\le (\log\log n)^{5/2+\varepsilon}$ for
sufficiently large $n$, and $\log \eta_{n/(\log\log n)^5}\sim (1/2)\log n$ as $n\rightarrow \infty$. Hence
 \bestar   \frac{l(\eta_n)-l(\eta_{n/(\log\log n)^5})}{l(\eta_n)}=o(1/\log\log n).~~~\Box \eestar

\begin{lemma} \label{de1} As $n \rightarrow \infty$, we have
\bestar \frac{\sum_{j=1}^n (|Z_j|+E|Z_j|)}{B_n/\sqrt{\log\log n}} \stackrel{P}{\rightarrow} 0.\eestar
\end{lemma}

\noindent {\bf Proof.}  Let $\tau_j=\eta_j(\log\log j)^3$ and $Z_j^{*}=X_jI(\eta_j<|X_j|<\tau_j)$. From
 the proof of Lemma 2 in Cs\"{o}rg\H{o} {\it et al.} (2003), we have $P(Z_j \ne Z_j^{*}, i.o.)=0$.
Hence, by Chebyshev's inequality, in order to prove Lemma \ref{de1}, we only need to
 prove that, as $n \rightarrow \infty$, \be && \sum_{j=1}^n E|Z_j^{*}|=o(B_n/\sqrt{\log\log n}),\label{hu1}\\
&& \sum_{j=1}^n EZ_j^{*2}=o(B_n^2/\log\log n), \label{hu2}\\
&& \sum_{j=1}^n E|X_j|I(|X_j|>\tau_j)=o(B_n/\sqrt{\log\log n}).\label{addhu1}
 \ee We only prove (\ref{hu1}) and (\ref{addhu1}), for the proof of (\ref{hu2}) is similar to that of (\ref{hu1}). Since $\eta_n$ is a
regularly varying function with index $1/2$, we have that for sufficiently large $n$,
$$\eta_{n/(\log\log n)^{16}}(\log\log n)^3\le \eta_{n/(\log\log n)^9}.$$
Also, similarly, by the fact that  $\sqrt{j}(\log\log j)^2/\sqrt{l(\eta_j)}$ is a regularly varying function with index $1/2$, we have that for
sufficiently large $n$,\bestar \max_{1\le j\le n/(\log\log n)^9}\frac{j}{\eta_j}=\max_{1\le j\le n/(\log\log n)^9}\frac{\sqrt{j}(\log\log
j)^2}{\sqrt{l(\eta_j)}}\le \frac{\sqrt{n}}{\sqrt{l(\eta_n)}(\log\log n)^2}. \eestar Hence, by using the same method as that in the proof of
Lemma \ref{lemma2}, we have \bestar && \sum_{j=1}^n E|Z_j^{*}|\le \sum_{j=1}^{n/(\log\log n)^{16}}
E|X_1|I(\eta_i<|X_1|<\eta_{n/(\log\log n)^9})\\
&&~~~~~~~~~~~~~~~~~~~ +n E|X_1|I(\eta_{n/(\log\log n)^{16}}<|X_1|<\eta_{n}(\log\log n)^3)\\
&&~~~~~~~~~~~~~~\le  \sum_{j=1}^{n/(\log\log n)^9} jE|X_1|I(\eta_j<|X_1|<\eta_{j+1})\\
&&~~~~~~~~~~~~~~~~~~~ +\frac{n (l(\eta_{n}(\log\log n)^3)-l(\eta_{n/(\log\log n)^{16}}))}{\eta_{n/(\log\log n)^{16}}}\\
&&~~~~~~~~~~~~~~ =o(B_n/(\log\log n)),~~n \rightarrow \infty.\eestar Thus  (\ref{hu1}) is proved.

Next, we prove (\ref{addhu1}). By the fact that $E|X|I(|X|\ge x)=o(1)l(x)/x$ as $x\rightarrow \infty$,
 \bestar \sum_{j=1}^n
E|X_j|I(|X_j|>\tau_j)=o(1) \sum_{j=1}^n  \frac{l(\tau_j)}{\tau_j} \le o(1) l(\tau_n) \sum_{j=1}^n  \frac{1}{\tau_j}. \eestar Since $1/\tau_n$
is a regularly varying function with index $-1/2$, by Tauberian theorem (see, for instance, Theorem 5 in Feller (1971), page 447), we have
$\sum_{j=1}^n \frac{1}{\tau_j} \sim 2{n}/{\tau_n}$ as $n \rightarrow \infty$. Hence, as $n \rightarrow \infty$, \bestar \sum_{j=1}^n
E|X_j|I(|X_j|>\tau_j)=
 o(1) \frac{n l(\tau_n)}{\tau_n}=o(1)B_n/(\log\log n). \eestar

 Thus  (\ref{addhu1}) is proved and the proof of Lemma (\ref{de1}) is complete. ~~$\Box$

\begin{lemma} \label{lemma1}
For all $t \in \mathbb{R}$, we have  \be \lim\limits_{n \rightarrow \infty} P\Big(a(n)\max_{1\le k< n}S_{k,n}^*/B_{k,n}\le
t+b(n)\Big)=\exp(-e^{-t}),\label{1205}\ee and \be \lim\limits_{n \rightarrow \infty} P\Big(a(n)\max_{1\le k< n}|S_{k,n}^*|/B_{k,n}\le
t+b(n)\Big)=\exp(-2e^{-t}).\label{add1205}\ee
\end{lemma}

\noindent {\bf Proof.} We only prove (\ref{1205}), since the proof of (\ref{add1205}) is similar.
 Since $l(x^2)\le 2^{C_0}l(x)$, by (42) in
Cs\"{o}rg\H{o} {\it et al.} (2003), there exist two independent Wiener processes $W^{(1)}$ and $W^{(2)}$ such that, as $n\rightarrow \infty$,
\be S_n^*-W^{(1)}(B_n^{2})=o(B_n/\sqrt{\log\log n})~~a.s. \label{hu11}\ee and \be \tilde{S}_n^*-W^{(2)}(B_n^{2})=o(B_n/\sqrt{\log\log n})~~a.s.
\label{hu12}\ee
 Define  $K_n=\exp\{\log^{1/3} n\}$ and
 \bestar W(n,t)=\left\{\begin{array}{ll} n^{-1/2}(W^{(1)}(nt)-t(W^{(1)}(n/2)+W^{(2)}(n/2))), & 0\le t\le 1/2,\\
 n^{-1/2}(-W^{(2)}(n-nt)+(1-t)(W^{(1)}(n/2)+W^{(2)}(n/2))), & 1/2<t\le 1. \end{array}\right.\eestar
Computing its covariance function, one concludes that $W(n,t)$ is a Brownian bridge in $0\le t\le 1$ for each $n\ge 1$.  Now, as $n \rightarrow
\infty$, we have
 \be \sqrt{\log \log n}\max_{K_n\le k\le
n/2}\Big|\frac{S_{k,n}^*}{B_{k,n}}-\frac{B_n^2W(B_n^2,B_k^2/B_n^2)}{\sqrt{B_k^2(B_n^2-B_k^2)}}\Big|\stackrel{P}{\rightarrow} 0.
 \label{add1}
\ee

To prove (\ref{add1}), we notice that for $k\le n/2$, \bestar
S_{k,n}^*=\frac{n}{k(n-k)}\Big(S_k^*-\frac{k}{n}(\tilde{S}_{n-[n/2]}^*+S_{[n/2]}^*)\Big).\eestar Hence, for $k\le n/2$, \be
&&\Big|\frac{S_{k,n}^*}{B_{k,n}}-\frac{B_n^2W(B_n^2,B_k^2/B_n^2)}{\sqrt{B_k^2(B_n^2-B_k^2)}}\Big| \le
|W(B_n^2,B_k^2/B_n^2)|\Big|\frac{nB_n}{k(n-k)B_{k,n}}- \frac{B_n^2}{\sqrt{B_k^2(B_n^2-B_k^2)}}\Big|\nonumber\\
&&~~~~~~~~~~~~~~~~~~~~~~+\frac{nB_n}{k(n-k)B_{k,n}}\Big|\frac{k(n-k)}{nB_n}S_{k,n}^*-W(B_n^2,B_k^2/B_n^2)\Big|\nonumber\\
&&~~~~~~~~~~~~~~~~~~~~:=L_1(k,n)+L_2(k,n). \label{1201}\ee

First, we estimate $L_1(k,n)$.  We have
 \bestar \frac{k^2(n-k)^2B_{k,n}^2}{n^2B_n^2}-
\frac{B_k^2(B_n^2-B_k^2)}{B_n^4}=\Big(\frac{B_k^2}{B_n^2}-\frac{k}{n}\Big)^2-\frac{k^2(B_n^2-B_{[n/2]}^2-B_{n-[n/2]}^2)}{n^2B_n^2}.\eestar
 Note
that   $ (k/n)^{5/8} \le B_k/B_n\le (k/n)^{3/8}$ holds for all $K_n \le k\le n$ and sufficiently large $n$ by the fact that $B_n$ is a
regularly varying function with index $1/2$. Then
 \bestar  \max_{K_n\le k\le
n/(\log\log n)^5}\frac{B_n^3}{B_k^3}\Big(\frac{B_k^2}{B_n^2}-\frac{k}{n}\Big)^2&\le& 2\max_{K_n\le k\le n/(\log\log n)^5}\Big(\frac{B_k}{B_n} +
\frac{B_n^3k^2}{B_k^3n^2}\Big)\\
&\le& 4(\log\log n)^{-5/8}.\eestar Also, by Lemma \ref{lemma2}, \bestar  &&\max_{n/(\log\log n)^5<k\le
n/2}\frac{B_n^3}{B_k^3}\Big(\frac{B_k^2}{B_n^2}-\frac{k}{n}\Big)^2 \le
 \max_{n/(\log\log n)^5<k\le  n/2}\frac{k^2B_n^3}{n^2B_k^3} \frac{(l(\eta_n)-l(\eta_{n/(\log\log n)^5}))^2}{l(\eta_n)^2}\\
&&~~~~~~~~~~~~~~~~~~~~~~~~~~~~~~~~~~~~~~~~=o(1/\sqrt{\log\log n}), ~~n\rightarrow \infty.\eestar Hence, as $n\rightarrow \infty$, \be
\sqrt{\log\log n} \max_{K_n\le k\le n/2}\frac{B_n^3}{B_k^3}\Big(\frac{B_k^2}{B_n^2}-\frac{k}{n}\Big)^2 \rightarrow 0. \label{log2}\ee
 Again by Lemma \ref{lemma2}, \bestar
\sqrt{\log\log n}\max_{K_n\le k\le n/2}\frac{B_n^3}{B_k^3}\frac{k^2(B_n^2-B_{[n/2]}^2-B_{n-[n/2]}^2)}{n^2B_n^2}\le {\sqrt{\log\log
n}}\frac{l(\eta_n)-l(\eta_{[n/2]})}{l(\eta_n)}\rightarrow 0, \eestar as $n\rightarrow \infty$. Thus, \be \sqrt{\log\log n}\max_{K_n\le k\le
n/2}\frac{B_n^3}{B_k^3}\Big|\frac{k^2(n-k)^2B_{k,n}^2}{n^2B_n^2}-\frac{B_k^2(B_n^2-B_k^2)}{B_n^4}\Big|\rightarrow 0, ~~n\rightarrow
\infty.\label{log1}\ee This implies that  for large $n$ and  all $K_n\le k\le n/2$, \bestar
\Big|\frac{k^2(n-k)^2B_{k,n}^2}{n^2B_n^2}-\frac{B_k^2(B_n^2-B_k^2)}{B_n^4}\Big|
&\le& \frac{1}{4} \frac{B_k^3}{B_n^3} \le \frac{1}{4}\frac{B_k^2(B_n^2-B_k^2)}{B_n^4}\frac{B_n B_{[n/2]}}{B_n^2-B^2_{[n/2]}} \\
&\le& \frac{1}{2}\frac{B_k^2(B_n^2-B_k^2)}{B_n^4}.\eestar Hence, for large $n$ and  all $K_n\le k\le n/2$, \be
\frac{(1/2)B_n^2}{\sqrt{B_k^2(B_n^2-B_k^2)}} \le \frac{nB_n}{k(n-k)B_{k,n}} \le \frac{2B_n^2}{\sqrt{B_k^2(B_n^2-B_k^2)}}. \label{add2}\ee
Noting that $|{1}/{\sqrt{x}}-{1}/{\sqrt{y}}|\le |x-y|/ (x\sqrt{y})$ for all $x, y>0$, it follows from (\ref{log1}) and (\ref{add2}) that
 \be && \sqrt{\log\log n}\max_{K_n \le k\le n/2}\Big|\frac{nB_n}{k(n-k)B_{k,n}}-
 \frac{B_n^2}{\sqrt{B_k^2(B_n^2-B_k^2)}}\Big|\nonumber\\
 &\le&  \sqrt{2\log\log n}\max_{K_n \le k\le n/2}\Big|\frac{k^2(n-k)^2B_{k,n}^2}{n^2B_n^2}-\frac{B_k^2(B_n^2-B_k^2)}{B_n^4}\Big|
 \Big(\frac{B_k^2(B_n^2-B_k^2)}{B_n^4}\Big)^{-3/2}\nonumber\\
 &\le&  4\sqrt{\log\log n}\max_{K_n \le k\le n/2}\frac{B_n^3}{B_k^3}\Big|\frac{k^2(n-k)^2B_{k,n}^2}{n^2B_n^2}-\frac{B_k^2(B_n^2-B_k^2)}{B_n^4}\Big|
\rightarrow 0. \label{add3} \ee By properties of Brownian motion,  \bestar \max_{K_n \le k\le n/2} |W(B_n^2,B_k^2/B_n^2)|&\le&
2B_n^{-1}\sup_{0\le t\le B_n^2} |W^{(1)}(t)|+
B_n^{-1}|W^{(2)}(B_n^2/2)|\\
&\stackrel{d}{=}& 2\sup_{0\le t\le 1} |W^{(1)}(t)|+|W^{(2)}(1/2)|.\eestar This together with (\ref{add3}) yields \be \sqrt{\log\log n}\max_{K_n
\le k\le n/2} L_1(k,n) \stackrel{P}{\rightarrow} 0,~~ n\rightarrow \infty. \label{add5}\ee

Next, we estimate $L_2(k,n)$.  By (\ref{hu11}) and (\ref{hu12}),  \bestar &&\Big|\frac{k(n-k)}{nB_n}S_{k,n}^*-W(B_n^2,B_k^2/B_n^2)\Big|\le
\frac{k}{nB_n}|W^{(1)}(B_n^2/2)-W^{(1)}(B_{[n/2]}^2)|\nonumber\\
&&~~~~~~~~~~~~~~~~~~~~+\frac{k}{nB_n}|W^{(2)}(B_n^2/2)-W^{(2)}(B_{n-[n/2]}^2)|\nonumber\\
&&~~~~~~~~~~~~~~~~~~~~+\Big|\frac{k}{n}-\frac{B_k^2}{B_n^2}\Big|
\frac{|W^{(1)}(B_n^2/2)|+|W^{(2)}(B_n^2/2)|}{B_n}+\frac{o_k(1)B_k}{B_n\sqrt{\log\log k}},~\eestar where $o_k(1)\rightarrow 0$ as $k\rightarrow
\infty$. Similarly to the proof of (\ref{log2}), we have \bestar \sqrt{\log\log n} \max_{K_n \le k\le
n/2}\frac{B_n}{B_k}\Big|\frac{B_k^2}{B_n^2}-\frac{k}{n}\Big| \rightarrow 0,~~n\rightarrow \infty.  \eestar This, together with (\ref{add2}) and
the fact that
$$\frac{|W^{(1)}(B_n^2/2)|+|W^{(2)}(B_n^2/2)|}{B_n}\stackrel{d}{=} |W^{(1)}(1/2)|+|W^{(2)}(1/2)|,$$
as $n\rightarrow \infty$, yields \bestar \sqrt{\log \log n}\max_{K_n\le k\le
n/2}\frac{nB_n}{k(n-k)B_{k,n}}\Big|\frac{k}{n}-\frac{B_k^2}{B_n^2}\Big| \frac{|W^{(1)}(B_n^2/2)|+|W^{(2)}(B_n^2/2)|}{B_n}
\stackrel{P}{\rightarrow} 0. \eestar
 Similarly to the proof of Lemma \ref{lemma2}, we have
\bestar &&\frac{\sqrt{\log\log n}}{B_n}|W^{(1)}(B_n^2/2)-W^{(1)}(B_{[n/2]}^2)|
\stackrel{d}{=}\sqrt{\log\log n}\Big(\frac{B_n^2/2-B_{[n/2]}^2}{B_n^2}\Big)^{1/2}|W^{(1)}(1)|\\
&&~~~~~~~~~~~~~~~~~=\sqrt{\log\log n}
\Big(\frac{(n/2)l(\eta_n)-[n/2]l(\eta_{[n/2]})}{nl(\eta_n)}\Big)^{1/2}|W^{(1)}(1)|\stackrel{P}{\rightarrow}0,~~n \rightarrow \infty.\eestar
Hence, by (\ref{add2}), as $n \rightarrow \infty$,
 \bestar \sqrt{\log \log n}\max_{K_n\le k\le n/2}\frac{nB_n}{k(n-k)B_{k,n}}
 \frac{k}{nB_n}|W^{(1)}(B_n^2/2)-W^{(1)}(B_{[n/2]}^2)|\stackrel{P}{\rightarrow} 0. \eestar
Similarly, as $n \rightarrow \infty$, \bestar \sqrt{\log \log n}\max_{K_n\le k\le n/2}\frac{nB_n}{k(n-k)B_{k,n}}
\frac{k}{nB_n}|W^{(2)}(B_n^2/2)-W^{(2)}(B_{n-[n/2]}^2)| \stackrel{P}{\rightarrow} 0. \eestar Also, by (\ref{add2}), as $n \rightarrow \infty$,
 \bestar \sqrt{\log \log n}\max_{K_n\le k\le n/2}\frac{nB_n}{k(n-k)B_{k,n}} \frac{o_k(1)B_k}{B_n\sqrt{\log\log k}}
\stackrel{P}{\rightarrow} 0. \eestar
 Hence
 \be \sqrt{\log\log
n}\max_{K_n \le k\le n/2} L_2(k,n) \stackrel{P}{\rightarrow} 0,~~n \rightarrow \infty. \label{add4} \ee
 Now (\ref{add1}) follows from (\ref{1201}), (\ref{add5}) and
(\ref{add4}). Now,  similarly, as $n \rightarrow \infty$, \bestar \sqrt{\log \log n}\max_{n/2<k\le
n-K_n}\Big|\frac{S_{k,n}^*}{B_{k,n}}-\frac{B_n^2W(B_n^2,B_k^2/B_n^2)}{\sqrt{B_k^2(B_n^2-B_k^2)}}\Big|\stackrel{P}{\rightarrow} 0. \eestar
Hence, as $n \rightarrow \infty$, \bestar \sqrt{\log \log n}\Big|\max_{K_n\le k\le n-K_n}\frac{S_{k,n}^*}{B_{k,n}}-\sup_{K_n\le k \le
n-K_n}\frac{W(B_n^2,t)}{\sqrt{(B_k^2/B_n^2)(1-B_k^2/B_n^2)}}\Big|\stackrel{P}{\rightarrow} 0. \eestar

Next, we will show that, as $n \rightarrow \infty$,  \be \sqrt{\log \log n}\Big|\sup_{\frac{B_{K_n}^2}{B_n^2}\le t \le
\frac{B_{n-K_n}^2}{B_n^2}}\frac{W(B_n^2,t)}{\sqrt{t(1-t)}}-\sup_{K_n\le k \le
n-K_n}\frac{W(B_n^2,t)}{\sqrt{(B_k^2/B_n^2)(1-B_k^2/B_n^2)}}\Big|\stackrel{P}{\rightarrow} 0.  \label{add6}\ee
 Write $$\Delta_n=\inf_{K_n+1\le
k\le n-K_n} \frac{B_k^2-B_{k-1}^2}{B_n^2}=\frac{l(\eta_{K_n})}{B_n^2} $$ and recall that $W(B_n^2, t)$ is a Brownian bridge in $t\in [0,1]$ for
each $n\ge 1$. Hence, to prove (\ref{add6}), we only need to show that, as $n \rightarrow \infty$, \bestar  \sqrt{\log \log
n}\sup_{\frac{B_{K_n}^2}{B_n^2}\le t,s \le \frac{B_{n-K_n}^2}{B_n^2}}\sup_{|t-s|\le
\Delta_n}\Big|\frac{W(t)-tW(1)}{\sqrt{t(1-t)}}-\frac{W(s)-sW(1)}{\sqrt{s(1-s)}}\Big|\stackrel{P}{\rightarrow} 0, \eestar where $W(t)$ is a
standard Brownian motion. This follows from results on the increments of a Brownian motion (see for instance  Cs\"{o}rg\H{o} and R\'{e}v\'{e}sz
(1981), Theorem 1.2.1) and by some basic calculations. We omit the details here.
 Hence, as $n \rightarrow \infty$, \be \sqrt{\log \log n}\Big|\max_{K_n\le k\le
n-K_n}\frac{S_{k,n}^*}{B_{k,n}}-\sup_{\frac{B_{K_n}^2}{B_n^2}\le t \le
\frac{B_{n-K_n}^2}{B_n^2}}\frac{W(B_n^2,t)}{\sqrt{t(1-t)}}\Big|\stackrel{P}{\rightarrow} 0. \label{add7} \ee
 By using (A.4.30) and (A.4.31) in
Cs\"{o}rg\H{o} and Horv\'{a}th (1997), as $n \rightarrow \infty$, we conclude \bestar &&(2\log \log B_n^2)^{-1/2} \sup_{1/B_n^2\le t \le
c(B_n^2)}\frac{W(B_n^2,t)}{\sqrt{t(1-t)}}
 \stackrel{P}{\rightarrow} \sqrt{5/12},\\
&&(2\log \log B_n^2)^{-1/2} \sup_{1-c(B_n^2)\le t \le 1/B_n^2 }\frac{W(B_n^2,t)}{\sqrt{t(1-t)}}
 \stackrel{P}{\rightarrow} \sqrt{5/12},
\eestar where $c(B_n^2)=\exp\{(\log B_n^2)^{5/12}\}/B_n^2$. Notice that $B_{K_n}^2/B_n^2 \le c(B_n^2)$ and $ B_{n-K_n}^2/B_n^2\ge 1-c(B_n^2)$
for sufficiently large $n$. Hence, as $n \rightarrow \infty$,  \be && a(B_n^2)\sup_{1/B_n^2\le t \le B_{K_n}^2/B_n^2}
\frac{W(B_n^2,t)}{\sqrt{t(1-t)}}-b(B_n^2)\stackrel{P}{\rightarrow}-\infty, \\
&& a(B_n^2)\sup_{ B_{n-K_n}^2/B_n^2\le t \le 1-1/B_n^2} \frac{W(B_n^2,t)}{\sqrt{t(1-t)}}-b(B_n^2)\stackrel{P}{\rightarrow}-\infty.\ee By
(A.4.29) and Theorem A.3.1 in Cs\"{o}rg\H{o} and Horv\'{a}th (1997), we arrive at \be \lim\limits_{n \rightarrow \infty}
P\Big(a(B_n^2)\sup_{1/B_n^2 \le t \le 1-1/B_n^2} \frac{W(B_n^2,t)}{\sqrt{t(1-t)}}\le t+b(B_n^2)\Big)=\exp(-e^{-t}). \label{add8}\ee Now, from
(\ref{add7})--(\ref{add8}) it follows that  for all $t \in \mathbb{R}$, \bestar \lim\limits_{n \rightarrow \infty} P\Big(a(B_n^2)\max_{K_n \le
k\le n-K_n}S_{k,n}^*/B_{k,n}\le t+b(B_n^2)\Big)=\exp(-e^{-t}).\eestar
 This, together with (\ref{1204}) below, implies  that for all $t \in \mathbb{R}$,
 \bestar \lim\limits_{n \rightarrow \infty} P\Big(a(B_n^2)\max_{1\le k<n}S_{k,n}^*/B_{k,n}\le
t+b(B_n^2)\Big)=\exp(-e^{-t}).\eestar Since, as $n\rightarrow \infty$, $\log \log B_n^2=\log\log n+o(1)$, we have \bestar &&a(n)\max_{1\le k<
n}S_{k,n}^*/B_{k,n}-b(n)\\
&&~~~~~~=\frac{a(n)}{a(B_n^2)}\Big(a(B_n^2)\max_{1\le k< n}S_{k,n}^*/B_{k,n}-b(B_n^2)\Big)+\frac{a(n)}{a(B_n^2)}b(B_n^2)-b(n)\\
&&~~~~~~=(1+o(1))\Big(a(B_n^2)\max_{1\le k<n}S_{k,n}^*/B_{k,n}-b(B_n^2)\Big)+o(1), \eestar which implies (\ref{1205}).  Lemma \ref{lemma1} is
proved. ~~$\Box$

\vskip 0.5cm

\noindent {\bf Proof of Theorem \ref{th1}}.  Write $K_n=\exp\{\log^{1/3} n\}$, and put \bestar && \Omega_1=\Big\{K_n< k \le n/4: \sum_{i=1}^k |Z_i|\le B_k/\log\log k\Big\},\\
&& \Omega_2=\Big\{K_n< k \le n/4: \sum_{i=1}^k |\tilde{Z}_i|\le B_k/\log\log k\Big\}.\eestar Define $\Omega'=\Omega_1\cup \{k: n/4< k\le
n/2\}$,~$\Omega''=\{k: n-k \in \Omega_2\}\cup \{k: n/2< k< 3n/4\}$ and
 $\Omega_1'=\{k: 2\le k\le n/4\}-\Omega_1,~\Omega_2'=\{k: 3n/4\le k\le n-2\}-\{k: n-k \in \Omega_2\}$.

Notice that, as $n\rightarrow \infty$, $S_{[nt]}/b_n\stackrel{d}{\rightarrow} W(t)$ and $V_n^2/b_n^2\stackrel{P}{\rightarrow} 1$, where $W$ is
a Brownian motion and $b_n$ is a regularly varying function with index $1/2$. Hence \bestar && \frac{\min_{k\le
n/4}(\tilde{V}^2_{n-[n/2]}+V^2_{[n/2]}-V_k^2-(\tilde{S}_{n-[n/2]}+S_{[n/2]}-S_k)^2/(n-k))}{b_n^2}\\&&~~~~~~~
 \ge \frac{\tilde{V}^2_{n-[n/2]}}{b_n^2}-\frac{3\tilde{S}_{n-[n/2]}^2+6(\max_{1\le k\le n/2}|S_k|)^2}
 {(n/2)b_n^2}\stackrel{P}{\rightarrow} 1/2,~~n\rightarrow \infty.\eestar
Notice that by the  self-normalized LIL of Griffin and Kuelbs (1989), as $n\rightarrow \infty$, we have \bestar \limsup_{n\rightarrow \infty}
\frac{|S_n|}{\sqrt{2\log\log n (V_n^2-S_n^2/n)}}=1~~a.s.\eestar Consequently,
 \bestar \frac{1}{\sqrt{2\log\log n}}\max_{2\le k\le K_n} \frac{|S_k|}{\sqrt{(V_k^2-S_k^2/k)}}\le
 \frac{\sqrt{2\log\log K_n}}{\sqrt{2\log\log n}}(1+o(1))=\sqrt{1/3}+o(1)~~a.s. \eestar
Similarly, by (18) in Cs\"{o}rg\H{o} {\it et al.} (2003), we conclude \bestar \frac{1}{\sqrt{2\log\log n}}\max_{ k>K_n~\mbox{and} ~k\in
\Omega_1'} \frac{|S_k|}{\sqrt{(V_k^2-S_k^2/k)}}\le \sqrt{1/2}+o(1)~~a.s.,~~n\rightarrow \infty. \eestar Thus, by noting that
$\frac{a+b}{\sqrt{c+d}}\le \frac{a}{\sqrt{c}}+\frac{b}{\sqrt{d}}$ holds for all $a,b,c,d>0$,
 \bestar && \frac{1}{\sqrt{2\log \log n}}\max_{k \in \Omega_1'} \frac{|S_{k,n}|}{\bar{V}_{k,n}}\le \frac{1}{\sqrt{2\log \log n}}\max_{k \in
\Omega_1'} \frac{n}{n-k}\frac{|S_{k}|}{\sqrt{V_k^2-S_k^2/k}}\\&&~~~~~~~+\frac{(|S_{[n/2]}|+|\tilde{S}_{n-[n/2]}|)/(b_n\sqrt{2\log \log n})}
{\min \limits_{k\le n/4}\sqrt{\tilde{V}^2_{n-[n/2]}+V^2_{[n/2]}-V_k^2-(\tilde{S}_{n-[n/2]}+S_{[n/2]}-S_k)^2/(n-k)}/b_n}\\&&~~~~~~~ \le
{2\sqrt{2}}/{3}+o_P(1),~~n\rightarrow \infty.\eestar This, as $n\rightarrow \infty$, implies \be a(n)\max_{k \in \Omega_1'}
\frac{|S_{k,n}|}{\bar{V}_{k,n}}-b(n)\stackrel{P}{\rightarrow}-\infty, \label{hu1201} \ee and, similarly \be a(n)\max_{k \in \Omega_2'}
\frac{|S_{k,n}|}{\bar{V}_{k,n}}-b(n)\stackrel{P}{\rightarrow}-\infty. \label{hu1202}\ee

Furthermore,  similarly, by using  (20) in Cs\"{o}rg\H{o} {\it et al.} (2003),  and by the facts that, as $n\rightarrow \infty$,
$S_n^{*}/B_n\stackrel{d}{\rightarrow} N(0,1)$ and $\limsup_{n\rightarrow \infty} S_n^{*}/ (2B_n^2\log\log n)^{1/2}=1~a.s.$ (by (\ref{hu11})),
 we infer \be a(n)\max_{k \in  \Omega_1'\cup\Omega_2'}
\frac{|S_{k,n}^{*}|}{B_{k,n}}-b(n)\stackrel{P}{\rightarrow}-\infty. \label{1204} \ee

Now, in order to prove Theorem \ref{th1}, we only need to show that, as $n\rightarrow \infty$,  \be a(n)\max_{k \in
\Omega'}\Big|\frac{S_{k,n}}{V_{k,n}}-\frac{S_{k,n}^*}{B_{k,n}}\Big|\stackrel{P}{\rightarrow} 0,\label{hu3} \ee and \be a(n)\max_{k \in
\Omega''}\Big|\frac{S_{k,n}}{V_{k,n}}-\frac{S_{k,n}^*}{B_{k,n}}\Big|\stackrel{P}{\rightarrow} 0.\label{hu4} \ee

In fact, if (\ref{hu3}) and (\ref{hu4}) hold true,  then it follows from (\ref{1204}) and Lemma \ref{lemma1} that, for all $t \in \mathbb{R}$,
\be \lim\limits_{n \rightarrow \infty} P\Big(a(n)\max_{k \in  \Omega'\cup\Omega''}{S_{k,n}}/{V_{k,n}}\le t+b(n)\Big)=\exp(-e^{-t}).
\label{addproof4} \ee And also by Lemma \ref{lemma1}, we obtain that \be \frac{1}{\sqrt{2\log\log n}}\max_{1\le
k<n}\frac{|S_{k,n}^*|}{B_{k,n}}\stackrel{P}{\rightarrow} 1, ~~~~n\rightarrow \infty. \label{addproof5}\ee

By noting that $$V_{k,n}^2\le \bar{V}_{k,n}^2\le \max\Big\{\frac{k}{k-1}, \frac{n-k}{n-k-1}\Big\}V_{k,n}^2,$$ and by applying (\ref{hu3}),
(\ref{hu4}) and (\ref{addproof5}), we get that,
 \be && a(n)\max_{k \in  \Omega'\cup\Omega''} \Big|\frac{S_{k,n}}{\bar{V}_{k,n}}-\frac{S_{k,n}}{V_{k,n}}\Big|
 \le \frac{a(n)}{\sqrt{K_n}}\max_{k \in  \Omega'\cup\Omega''}\frac{|S_{k,n}|}{V_{k,n}} \nonumber \\
 &&~~~~~~~~~~\le \frac{a(n)}{\sqrt{K_n}}\max_{k \in  \Omega'\cup\Omega''}\Big|\frac{S_{k,n}}{V_{k,n}}-\frac{S_{k,n}^*}{B_{k,n}}\Big|
 +\frac{a(n)}{\sqrt{K_n}}\max_{k \in  \Omega'\cup\Omega''}\frac{|S_{k,n}^*|}{B_{k,n}}\nonumber\\
 && ~~~~~~~~~~\stackrel{P}{\rightarrow} 0,~~n\rightarrow \infty.
  \label{addproof3} \ee
This, together with (\ref{hu1201}), (\ref{hu1202}) and (\ref{addproof4}), yields Theorem \ref{th1}.

Now we go to prove (\ref{hu3}) and (\ref{hu4}).
 We only prove (\ref{hu3}), since the proof of (\ref{hu4}) is similar. Clearly, we have \be &&
a(n)\max_{k \in \Omega'}\Big|\frac{S_{k,n}}{V_{k,n}}-\frac{S_{k,n}^*}{B_{k,n}}\Big| \le a(n)\max_{k\in \Omega'}
\Big|\frac{S_{k,n}}{V_{k,n}}-\frac{S_{k,n}}{B_{k,n}}\Big|+a(n)\max_{k\in
\Omega'}\Big|\frac{S_{k,n}-S_{k,n}^*}{B_{k,n}}\Big|\nonumber\\
&&~~~~~~~~~~~~~~~\le  a(n) \max_{k\in \Omega'} \Big|\frac{S_{k,n}}{V_{k,n}}\frac{V_{k,n}^2-B_{k,n}^2}{B^2_{k,n}}\Big| +a(n)\max_{k\in
\Omega'}\Big|\frac{S_{k,n}-S_{k,n}^*}{B_{k,n}}\Big|.\label{hu6}\ee By the self-normalized LIL of Griffin and Kuelbs (1989), we get that, as
$n\rightarrow \infty$,
$$\sup_{K_n\le k\le n/2}\frac{V_{k,n}^2} {{V_k^2}/{k^2}+ (\tilde{V}^2_{n-[n/2]}+V^2_{[n/2]}-V_k^2)/{(n-k)^2}} \rightarrow 1 ~~a.s.$$ Hence,
 for sufficiently large $n$,
 \be  a(n) \max_{k\in \Omega'}
\Big|\frac{S_{k,n}}{V_{k,n}}\frac{V_{k,n}^2-B_{k,n}^2}{B^2_{k,n}}\Big| &\le& 2a(n) \max_{k\in \Omega'}
\Big|\frac{S_{k}}{V_{k}}\frac{V_{k,n}^2-B_{k,n}^2}{B^2_{k,n}}\Big|\nonumber\\
&&+2a(n) \frac{V_{n}}{\tilde{V}_{n-[n/2]}}\max_{k\in \Omega'}\Big|\frac{S_{[n/2]}-S_k}{V_n}\frac{V_{k,n}^2-B_{k,n}^2}{B^2_{k,n}}\Big|
\nonumber\\
&&+2a(n) \max_{k\in \Omega'} \frac{|\tilde{S}_{n-[n/2]}|}{\tilde{V}_{n-[n/2]}}\Big|\frac{V_{k,n}^2-B_{k,n}^2}{B^2_{k,n}}\Big|.\label{hu7}\ee
 Since $EX=0$ and  $E|X_1|^r<\infty$ for any $1<r<2$,  it follows from the  Marcinkiewicz-Zygmund  strong law of large number (c.f. Chow and
Teicher (1978), page 125) that ${S_n}/{n^{1/r}}\rightarrow 0~a.s.$ Hence, as $n\rightarrow \infty$, $$\frac{(\log\log n)
S_n^2}{nB_n^2}\rightarrow 0 ~~a.s.$$ Note that for $n/4< k\le n/2$,
$$\frac{\sum_{j=1}^k (Z_j^2+|EZ_j|^2)/k^2}{B_{[n/2]}^2/(n-k)^2}\le 9\frac{\sum_{j=1}^k (Z_j^2+|EZ_j|^2)}{B_{[n/2]}^2},$$ and, by Lemma
\ref{de1}, \bestar && \frac{\sum_{j=1}^n |Z_j|^2}{B_n^2/\log\log n}\le \Big(\frac{\sum_{j=1}^n |Z_j|}{B_n/\sqrt{\log\log n}}\Big)^2
\stackrel{P}{\rightarrow} 0, \\
&&  \frac{\sum_{j=1}^n |EY_j|^2}{B_n^2/\log\log n}=\frac{\sum_{j=1}^n |EZ_j|^2}{B_n^2/\log\log n}\le \Big(\frac{\sum_{j=1}^n
|EZ_j|}{B_n/\sqrt{\log\log n}}\Big)^2 \rightarrow\, 0,~~n\rightarrow \infty. \eestar Now, by (40) of Cs\"{o}rg\H{o} {\it et al.} (2003), we
have
 \be &&(\log\log n) \max_{k\in \Omega'} \Big|\frac{V_{k,n}^2-B_{k,n}^2}{B^2_{k,n}}\Big|
\le 3\max_{k \in \Omega'} \frac{\log\log k|\sum_{j=1}^k (Y_j^2-EY_j^2)|}{B_k^2}\nonumber\\
&&~~~~~~~~~~~~~~~~+\frac{\log\log n|\sum_{j=1}^{[n/2]} (Y_j^2-EY_j^2)|}{B_{[n/2]}^2}+\frac{\log\log n|\sum_{j=1}^{n-[n/2]}
(\tilde{Y}_j^2-EY_j^2)|}{B_{n-[n/2]}^2}\nonumber\\
&&~~~~~~~~~~~~~~~~+3\max_{k \in \Omega_1} \frac{\log\log k\sum_{j=1}^k (Z_j^2+|EY_j|^2)}{B_k^2}+10\frac{\log\log n \sum_{j=1}^{[n/2]}
(Z_j^2+|EY_j|^2)}{B_{[n/2]}^2}\nonumber\\
&&~~~~~~~~~~~~~~~~+\frac{\log\log n\sum_{j=1}^{n-[n/2]} (\tilde{Z}_j^2+|EY_j|^2)}{B_{n-[n/2]}^2}+ 12\max_{k \in \Omega'} \frac{(\log\log
k)S_k^2}{kB_k^2}\nonumber\\
&&~~~~~~~~~~~~~~~~+3\frac{(\log\log n) \tilde{S}_{n-[n/2]}^2}{(n/2)B_{n-[n/2]}^2}+3\frac{(\log\log n) S_{[n/2]}^2}{(n/2)B_{[n/2]}^2}
 \stackrel{P}{\rightarrow} 0,~~ n\rightarrow \infty. \label{hu5}\ee
 By the self-normalized  LIL of Griffin and Kuelbs (1989), we conclude \be \max_{k\le n/2}
 \frac{|S_{[n/2]}-S_k|}{V_n\sqrt{2\log\log n}}\le \frac{2\max_{k\le
n/2}|S_k|}{V_n\sqrt{2\log\log n}}\le 2~~a.s,~~n\rightarrow \infty. \ee By the facts that $V_n^2/b_n^2\stackrel{P}{\rightarrow} 1$ and
$\tilde{V}_n^2/b_n^2\stackrel{P}{\rightarrow} 1$, as $n\rightarrow \infty$, we get \be \frac{V_{n}}{\tilde{V}_{n-[n/2]}}=
\frac{V_n}{b_n^2}\frac{b_{n-[n/2]}^2}{\tilde{V}_{n-[n/2]}}\frac{b_n^2}{b_{n-[n/2]}^2}\stackrel{P}{\rightarrow} 2. \label{add9}\ee Thus, by
using (\ref{hu7})-(\ref{add9}) and applying again the self-normalized  LIL of Griffin and Kuelbs (1989), as $n\rightarrow \infty$, we arrive at
 \be a(n) \max_{k\in \Omega'}
\Big|\frac{S_{k,n}}{V_{k,n}}\frac{V_{k,n}^2-B_{k,n}^2}{B^2_{k,n}}\Big|\stackrel{P}{\rightarrow} 0.\label{hu8}\ee

 Similarly to the proof of (\ref{hu5}),
by using  Lemma \ref{de1}, we have \be && a(n)\max_{k\in \Omega'}\Big|\frac{S_{k,n}-S_{k,n}^*}{B_{k,n}}\Big|
\le \sqrt{3}\max_{k \in \Omega_1} \frac{\sqrt{\log\log k}\sum_{j=1}^k (|Z_j|+|EZ_j|)}{B_k}\nonumber\\
&&~~~~~~~~~~~~~+4\frac{\sqrt{\log\log n} \sum_{j=1}^{[n/2]} (|Z_j|+|EZ_j|)}{B_{[n/2]}}+\frac{\sqrt{\log\log n}\sum_{j=1}^{n-[n/2]}
(|\tilde{Z}_j|+|EZ_j|)}{B_{n-[n/2]}}\nonumber \\
&&~~~~~~~~~~~~\stackrel{P}{\rightarrow} 0,~~n\rightarrow \infty.\label{hu9}\ee
 Now  (\ref{hu3}) follows from  (\ref{hu6}), (\ref{hu8}) and  (\ref{hu9}). This also completes the proof of Theorem \ref{th1}. ~~$\Box$

\end{document}